\documentclass{article}
\usepackage{graphicx}
 \usepackage{mathptmx}
\usepackage{amsmath, amstext,amssymb,amsfonts}
\usepackage[english]{babel}
\usepackage{delarray}
\usepackage{mathptmx}
\usepackage{amsmath, amstext,amssymb,amsfonts}
\usepackage[english]{babel}
\usepackage{delarray}

\DeclareFontFamily{U}{mathx}{\hyphenchar\font45}
\DeclareFontShape{U}{mathx}{m}{n}{
      <5> <6> <7> <8> <9> <10>
      <10.95> <12> <14.4> <17.28> <20.74> <24.88>
      mathx10
      }{}
\DeclareSymbolFont{mathx}{U}{mathx}{m}{n}
\DeclareFontSubstitution{U}{mathx}{m}{n}
\DeclareMathAccent{\widecheck}{0}{mathx}{"71}
\DeclareMathAccent{\wideparen}{0}{mathx}{"75}

\numberwithin{equation}{section}

\newcommand{\Co}{\mathbb C} \newcommand{\R}{\mathbb R}   \newcommand{\Z}{\mathbb Z}
       
\newcommand{\T}{\mathbb{T}}            \newcommand{\be}{B_{\sigma}^{p}}  
\newcommand{\bi}{B^1_{\sigma}}

\date{}
\title{Approximation of entire functions of exponential type by trigonometric sums}
\author{\large{ Saulius  Norvidas}}
\date{\footnotesize Institute of Data Science and Digital Technologies, Vilnius University, \\ Akademijos str. 4, Vilnius LT-04812, Lithuania\\
 ({\rm{e-mail: norvidas{@}gmail.com}})}
\begin{document}

\maketitle
 {{ {\bf Abstract}}}
 Let $\sigma>0$. For $1\le p\le \infty$,  the Bernstein space $\be$ is a  Banach space  of all  $f\in L^p(\R)$ such that $f$ is bandlimited to $\sigma$; that  is,  the distributional Fourier transform of $f$ is supported in $[-\sigma, \sigma]$.   We study   the   approximation of\ $f\in \be$\ by finite trigonometric  sums
\[
P_{\tau}(x)=\chi_{\tau}(x)\cdot  \sum_{|k|\le \sigma\tau/\pi}c_{k,\tau} e^{i\frac{\pi}{\tau}k x }
\]
in  $L^p$ norm on $\R$ as\ $\tau\to\infty$,\ where\ $\chi_{\tau}$ denotes  the indicator function of  $[-\tau, \tau]$.



\section{ Introduction }

{\large{\ \ \ \ \  The  Fourier transform  of   $f\in L^1(\R)$  is defined here by
\[
\hat{f}(x)=\int_{\R} f(t)e^{-ix t} dt.
\]
We normalize  the  inverse Fourier transform
\[
\check{f}(\xi)=\frac1{2\pi} \int_{\R} f(t)e^{i\xi t} dt
\]
so that the inversion formula $\widehat{(\check{f})}=f$  holds for suitable  $f$. If  $f\in L^p(\R)$ and $p>1$, then we understood  the Fourier transform    $\hat{f}$ in a distributional sense.
  Let $\sigma>0$. A function  $f$ is called bandlimited to $\sigma$ ( to $[-\sigma, \sigma]$) if
$\hat{f}$\ vanishes outside $[-\sigma, \sigma]$. We denote by $\be$ the Bernstein space of all functions $f\in L^p(\R)$ such that $f$ is bandlimited to $\sigma$. The space  $\be$ is equipped
with the norm}

 \[
\|f\|_p=\biggl(\int_{-\infty}^{\infty}|f(x)|^p dx\biggr)^{1/p}, \quad {\text{if}} \ 1\le p<\infty
\]
and
\[
\|f\|_{\infty}= {\mbox{ess\  sup}}_{-\infty<x<\infty} |f(x)|, \quad {\text{if}} \  p=\infty.
\]

\large{
By the Paley--Wiener--Schwartz theorem, any $f\in\be$ has an extension  onto the complex plane $\Co$ to an entire function of exponential type at most $\sigma$. That
is, for every\ $\varepsilon >0$\ there  exists \ $M>0$\ such that
\[
|f(z)|\leqslant  M e^{(\sigma+\varepsilon)|z|}
\]
for all\ $ z\in\Co$.\    There are several characterizations of  $\be$. For example, a function $f$ belongs to $\be$ if and only if $f$ is an entire function satisfying the Plancherel--Polya inequality
\begin{equation}
\biggl(\int_{-\infty}^{\infty}|f(x+iy)|^p\ dx\biggr)^{1/p}\le \|f\|_p  e^{\sigma |y|},
\end{equation}
where $y \in\R$. Note   that $1\le p\le r\le\infty$ implies (see [3, p. 49, Lemma 6.6])
\begin{equation}
\be \subset B_{\sigma}^r.
\end{equation}
\ \ \ \ \ Given\ $f\in B_{\sigma}^2$,\    Martin\ [9]\ studied  the  approximation of $f$ in a suitable fashion by   trigonometric polynomials.
To be precise, for any\ $\tau>0$,\ set
\begin{equation}
f_{\tau}(x)=\sum_{k=-N}^{N}c_{k,\tau} e^{i\frac{\pi}{\tau}kx},
\end{equation}
where  \ $c_{k,\tau}$\ is the Fourier coefficient of\ $f$\
\[
c_{k,\tau}=\frac1{2\tau}\int_{-\tau}^{\tau}f(t)e^{-i\frac{\pi}{\tau}k t}\,dt
\]
and $N=[\sigma \tau/\pi]_{\Z}.$  Here,   $[\omega]_{\Z}$ denotes    the integer part of a number $\omega$.
Define the  $\tau$-truncated version of   (1.3) by
\begin{equation}
\varphi_{f,\tau}:= f_{\tau}\cdot \chi_{[-\tau, \tau]},
\end{equation}
where $\chi_{[-\tau, \tau]}$\   is  the indicator (characteristic) function of  $[-\tau, \tau]$.   Recall that a function $f$
 is said to be strictly bandlimited to $[-\sigma,\sigma]$ if it is bandlimited
to $[-\varrho, \varrho]$ with $0\le \varrho<\sigma$.  It was proved  in [9] that if $f\in B_{\sigma}^2$ is strictly bandlimited to $[-\sigma,\sigma]$, then    (1.4) converges to\ $f$\
as\ $\tau\to\infty$ in\ $L^2$ norm on\ $\mathbb{R}$\ and on  horizontal lines in\ $\Co$.   This result was  extended  in [12] to
 the Bernstein spaces  $B^{2}_K$ consisting of  functions of several variables bandlimited to a compact set $K\subset\R^n$.

\ \ \ \ \   We now state the main result of this paper.

\ \ \ \ \ { THEOREM 1.}\quad\textit{Let $f\in\be$. If \ $1< p<\infty$, then }
\begin{equation}
\lim_{\tau\to\infty} \|f-\varphi_{f,\tau}\|_{L^p(\R)}=0.
\end{equation}

\ \ \ \  It is easily be checked  that (1.5)  is no longer true  for $p=\infty$ in general. Indeed, if $\sigma=1$ and
$f(x)=e^{ix}$, then
\[
f_{\tau}(x)=\sum_{k=-N}^{N}(-1)^k \frac{\sin \tau}{\tau-\pi k}  e^{i\frac{\pi}{\tau}kx}.
\]
Thus, if $\tau_{m}=\pi/2+2\pi m$, $m\in \Z$, then
\[
 \|f-  f_{\tau_m}\|_{L^{\infty}[-\tau_m,\tau_m]}\ge \Bigl|(f- f_{\tau_m})(\tau_m)\Bigr|\ge \Im (f- f_{\tau_m})(\tau_m)=1=\|f\|_{L^{\infty}(\R)}.
\]

\ \ \ \ \  Approximation by trigonometric polynomials in $B_{\sigma}^{\infty}$ was also considered.  Recall (see [2, p. 549]) that a sequence $(\psi_m)_{1}^{\infty}$, $\psi_m\in L^{\infty}(\R)$, is said   to converges narrowly to  $\psi\in L^{\infty}(\R)$ if $(\psi_m)_{1}^{\infty}$ converges uniformly to  $\psi$ on compact sets of $\mathbb{R}$ and  $\lim_{m\to\infty}\|\psi_m\|_{L^{\infty}(\R)}=\|\psi\|_{L^{\infty}(\R)}$. Lewitan [8] was the first to define   the following polynomials
\begin{equation}
\tilde{f}_{\tau}(x):= \sum_{k\in \mathbb{Z}}f(x+k\tau)\frac{\sin^2(\frac{x}{\tau}+k)}{(\frac{x}{\tau}+k)^2}, \ \tau>0.
\end{equation}
It was proved  in [8] that if $f\in B^{\infty}_{ \sigma}$, then  $\tilde{f}_{\tau}$ converges narrowly to $f$.  Krein [6] proposed  a new proof of this statement,
calculated  the degree of  approximation, and given several applications of (1.6). H\"{o}rmander [4]  considered more general trigonometric polynomials of such type  (for more  detailed history and references see [1, p.p. 146--152]). Schmeisser [14]   proved  similar convergence theorems for $\tau$-truncated version of the  Lewitan ( Lewitan-Krein-H\"{o}rmander )  polynomials (1.6)  with respect to $L^p$ norms on lines in $\mathbb{C}$ parallel to $\R$   and  with respect to $l^p$ norms in the case when  $f\in B^p_{\sigma}$,  $p\in[1, \infty)$.  The Lewitan  polynomials  of several variables   were also studied (see [11] for more details and the literature cited there).

\ \ \ \  \ Several well-known theorems were proved by using  (1.6). The  proofs are based  usual on suitable  properties of  trigonometric polynomials and  on the  narrowly convergence of (1.6).  For example, Krein [6]              extended the  Fejer and Riesz theorem on non-negative trigonometric polynomials to any non-negative  entire functions of exponential type $\sigma$: such an entire function $G$ can be written  on $\R$ as $|g|^{2}$, where $g$ is  some entire function  of exponential type at most $\sigma/2$. Other examples  can be found in  [4-5], [13]  and [15, p. 244].

\ \ \ \ It turns out that  (1.4)  as well as the Lewitan polynomials (1.6)  converge  narrowly. Moreover, the following  stronger statement  holds.

\ \ \ \  { THEOREM 2.}\quad\textit{Let $f\in\be$ and \ $1< p<\infty$. Then }
\begin{equation}
\lim_{\tau\to\infty} \|f-\varphi_{f,\tau}\|_{L^{\infty}(\R)}=0.
\end{equation}

{\centering
\section{ Preliminaries and proofs }}

\ \ \ \ \ Let $\Z$ denote the group of integers and  $\mathbb{T}=(-\pi, \pi]=\R/(2\pi\Z)$. We  write (1.3) in the form
\begin{equation}
f_{\tau}(x)=\frac1{2\tau}\int_{-\tau}^{\tau}f(t)D_N\Bigl(\frac{\pi}{\tau}(x-t)\Bigr) dt,
\end{equation}
where $D_N$ is the Dirichlet kernel defined by
\[
D_N(\xi)=\sum_{k=-N}^N e^{i k\xi}.
\]
Therefore,  setting
\begin{equation}
A_{\tau}(f)=\varphi_{f,\tau},
\end{equation}
we may assume that (1.4) defines on $\be$  the  one-parametric  family $(A_{\tau})_{\tau>0}$ of   bounded linear operators $A_{\tau}: \be\to L^p(\R)$.

\ \ \ \ \ { LEMMA 1.}\quad\textit{Let $p\in(1,\infty)$. Then there exists a constant $a(p)>0$ such that
\begin{equation}
\| A_{\tau}(f)\|_{L^p(\R)}\le a(p) \|f\|_p
\end{equation}
for every $f\in\be$ and all $\tau>0$.}

\ \ \ \ \ { PROOF.}\quad Let us define  $u_f(t)= f(\tau t/\pi )$, $t\in (-\pi, \pi]$. We will consider $u_f$ as an element of $L^p(\T)$. It follows from (2.1) and (2.2) that
\begin{gather}
\| A_{\tau}(f)\|^p_{L^p(\R)}=\| A_{\tau}(f)\|^p_{L^p[-\tau,\tau]}=\frac{\tau}{2^p\pi^{1+p}}\int_{-\pi}^{\pi}\Bigl|\int_{-\pi}^{\pi}f\Bigl(\frac{\tau}{\pi}\xi\Bigr)
D_N(y-\xi)\,d \xi\Bigr|^{p} dy\nonumber \\
=\frac{\tau}{2^p\pi^{1+p}}\int_{-\pi}^{\pi}\Bigl|\int_{-\pi}^{\pi}u_f(\xi)
D_N(y-\xi)\,d \xi\Bigr|^{p} dy=\frac{\tau}{2^p\pi^{1+p}}\|u_f\ast D_N\|^p_{L^p(\mathbb{T})},
\end{gather}
where  $N=[\sigma\tau/\pi]_{\mathbb{Z}}$. It is well known (see for example  [16,  p. 26, theorem 2.1.3]) that if $1<p<\infty$, then  there exists such $c(p)>0$ that $\|\varphi\ast D_N\|_{L^p(\T)}\le c(p)\|\varphi\|_{L^p(\T)}$ for every $\varphi\in L^p(\T)$ and all $N=0,1,2,\dots$. If we combine this with (2.4), we get
\[
\| A_{\tau}(f)\|_{L^p(\R)}\le \frac{c(p)}{2\pi}\ \|f\cdot\chi_{[-\tau,\tau]}\|_{L^p(\R)}
\]
 for all $\sigma, \tau>0$ and every $f\in \be$. Thus we obtain   (2.3)  with $a(p)=c(p)/(2\pi)$.

\ \ \ \ \ Next we estimate  the difference  $ \sin(ax)/x -\alpha D_N(bx)$ for some  $a, b$ and $\alpha$.

\ \ \ \ \ { LEMMA 2.}\quad\textit{Let $\tau, \sigma>0$ and  let $N=[\sigma\tau/\pi]_{\Z}$.  If  $0\le \delta<1$, then }
\begin{equation}
\max_{|v|\le (1+\delta)\tau}\biggl|\frac{\sin \sigma v}{\pi v}- \frac1{2\tau}D_N\Bigl(\frac{\pi}{\tau} v\Bigr)\biggr|< \frac{3+\omega\Bigl(\frac{\pi}2(1+\delta)\Bigr)}{2\tau},
\end{equation}
where
\begin{equation}
\omega(t)=\frac1{t}- \cot t.
\end{equation}

\ \ \ \ \ { PROOF.}\quad Substituting $2\tau y/\pi$ for $v$ in
\[
\frac{\sin \sigma v}{\pi v}- \frac1{2\tau}D_N\Bigl(\frac{\pi}{\tau} v\Bigr) ,
\]
we get
\[
\frac{\sin \sigma v}{\pi v}-\frac1{2\tau}\sum_{k=-N}^N e^{i\frac{\pi}{\tau} k v}=\frac1{2\tau}\biggl(\frac{\sin\bigl(\frac{2\sigma\tau}{\pi}y\bigl)}{y}-
\frac{\sin\bigl((2N+1)y\bigl)}{\sin y} \Biggl)
\]
\[
=\frac1{2\tau}\biggl(\frac{\sin \frac{2\sigma\tau}{\pi} y -\sin 2Ny}{y} -\cos 2Ny +\omega(y)\sin 2Ny\biggl),
\]
\[
=\frac1{2\tau}\biggl(\frac{2\sin \Bigl( \Bigl(\frac{\sigma\tau}{\pi}- \Bigl[\frac{\sigma\tau}{\pi}\Bigr]_{\Z}\Bigr) y \Bigl)\cos  \Bigl( \Bigl(\frac{\sigma\tau}{\pi}+ \Bigl[\frac{\sigma\tau}{\pi}\Bigr]_{\Z}\Bigr) y \Bigl)}{y} -\cos 2Ny +\omega(y)\sin 2Ny\biggl),
\]
where $\omega$ is defined in (2.6).
If $|v|\le (1+\delta)\tau$, then  $|y|\le \pi(1+\delta)/2<\pi$. Thus, to conclude the proof, it remains to note that   the function (2.6) is   odd on $(-\pi,\pi)$ and  increases on $(0,\pi)$.

\ \ \ \ \ { PROOF OF THEOREM 1.}\quad First, we claim that if  $1\le p<\infty$, then  $B_{\sigma}^1$ is a dense subset of  $\be$. Indeed, let $f\in\be$, $1\le p<\infty$. If  $0<\varrho<1$, then the function
\[
f_{\varrho}(x)= \frac{\sin^2 \varrho x}{(\varrho x)^2}f\bigl((1-\varrho^2) x\bigr)
\]
belongs both to $\be$ and to $B_{\sigma}^1$. Moreover, as $\varrho\to 0$, these $f_{\varrho}$ tend to $f$ in $L^p$ norm on $\R$. Finally, by (1.2), our claim is proved.

\ \ \ \ \ \ Now, by the uniform boundedness principle (the Banach--Steinhaus theorem) and   Lemma 1,  it suffices to check  (1.5) on any function in $ B_{\sigma}^1$.  Let  $f\in  B_{\sigma}^1$. Fix $0<\delta<1$.  In order to prove (1.5), we estimate the $L^p(-\tau, \tau)$-norm of the following  functions:
\begin{equation}
F_1(x)= f(x)- \int_{-\delta\tau}^{\delta\tau} f(t)\frac{\sin \sigma(x-t)}{\pi(x-t)} dt,
\end{equation}
\begin{equation}
F_2(x)=  \int_{-\delta\tau}^{\delta\tau} f(t)\biggl(\frac{\sin \sigma(x-t)}{\pi(x-t)} - \frac1{2\tau}D_N\Bigl(\frac{\pi}{\tau}(x-t)\Bigr)\biggr) dt,
\end{equation}
and
\begin{equation}
F_3(x)=  \frac1{2\tau}\int_{\delta\tau\le t\le \tau} f(t)D_N\Bigl(\frac{\pi}{\tau}(x-t)\Bigr)\biggr) dt,
\end{equation}
where   $N=[\sigma\tau/\pi]_{\mathbb{Z}}$.

\ \ \ \ Let us  begin with  $F_2$.  By Minkowski's inequality
\begin{equation}
\Bigl(\int_a^b\Bigl|\int_c^d u(x,y)\,dy\Bigr|^pdx\Bigr)^{1/p}\le \int_c^d \Bigl(\int_a^b| u(x,y)|^p\,dx\Bigr)^{1/p}dy
\end{equation}
 (see [15, p. 592, inequality (12)]), we have
\[
\Bigl(\int_{-\tau}^{\tau} |F_2(x)|^p dx\Bigr)^{\frac{1}{p}}\le \int_{-\delta\tau}^{\delta\tau}\biggl(\int_{-\tau}^{\tau}\Bigl|\frac{\sin \sigma(x-t)}{\pi(x-t)}-\frac1{2\tau}D_N\Bigl(\frac{\pi}{\tau}(x-t)\Bigr)\Bigr|^p dx\biggr)^{\frac1{p}} |f(t)| dt.
\]
Since $|x-t|\le \tau(1+\delta)$, it follows from (2.5) that
\begin{equation}
\Bigl(\int_{-\tau}^{\tau} |F_2(x)|^p dx\Bigr)^{\frac{1}{p}}\le \frac{3+\omega\bigl(\frac{\pi}2(1+\delta)\bigr)}{(2\tau)^{1-1/p}}\int_{-\delta\tau}^{\delta\tau} |f(t)|\,dt.
\end{equation}

\ \ \ \ \ Now we consider  $F_1$. The function $\sin \sigma x/\pi x$ belongs to  each $B_{\sigma}^q$ with $1<q\le\infty$. In addition, we have
$\hat{f}=\chi_{[-\sigma,\sigma]}$. Therefore, if $f\in\be$,  where $1\le p<\infty$, then
\[
f(x)=\int_{-\infty}^{\infty} f(t)\frac{\sin \sigma(x-t)}{\pi(x-t)}\, dt
\]
for   $  x\in\R$  (see  [3, p. 50, theorem 6.11]). Applying this  representation to $f$  in (2.7) and using   (2.10),  we have
\[
\Bigl(\int_{-\tau}^{\tau} |F_1(x)|^p dx\Bigr)^{\frac{1}{p}}\le \int_{|t|\ge \delta\tau}\biggl(\int_{-\tau}^{\tau}\Bigl|\frac{\sin\sigma(x-t)}{\pi(x-t)}\Bigr|^p dx\biggr)^{\frac1p}|f(t)|\,dt
\]
\begin{equation}
\le \Bigl\|\frac{\sin \sigma v}{\pi v}\Bigr\|_{L^p(\R)} \int_{|t|\ge \delta\tau} |f(t)|\, dt.
\end{equation}
A similar argument shows that
\begin{equation}
\Bigl(\int_{-\tau}^{\tau} |F_3(x)|^p dx\Bigr)^{\frac{1}{p}}\le \frac1{2\tau}\Bigl(\int_{-\tau}^{\tau}\Bigl|D_N\Bigl(\frac{\pi}{\tau}\zeta\Bigr)\Bigr|^p\, d\zeta\Bigr)^{\frac1{p}}\cdot
\int_{\delta\tau\le t\le \tau}|f(t)|\,dt.
\end{equation}
It is known (see, for example, [10, p. 153]) that if $1<p<\infty$, then  there exists $0<\alpha(p)<\infty$ such that
 \[
 \|D_N\|_{L^p(\T)}\le \alpha(p)(2N+1)^{1-1/p}.
 \]
Therefore, since  $N\le \sigma\tau/\pi$, it follows from (2.13) that
\begin{equation}
\Bigl(\int_{-\tau}^{\tau} |F_3(x)|^p dx\Bigr)^{\frac{1}{p}}\le\alpha(p)\Bigl(\frac{2\sigma}{\pi}+\frac1{\tau}\Bigr)^{1-1/p}\int_{\delta\tau\le t\le \tau}|f(t)|\,dt.
\end{equation}
We have
\[
F_1+F_2-F_3=f -f_{\tau}.
\]
Hence, combining our estimates (2.11), (2.12), and (2.14) we obtain
\[
\lim_{\tau\to\infty}\|f\cdot\chi_{[-\tau,\tau]}-\varphi_{\tau, f}\|_{L^p(\R)}=\lim_{\tau\to\infty}\|(f-f_{\tau})\chi_{[-\tau,\tau]}\|_{L^p(\R)}=0.
\]
Finally,  the triangle inequality gives us (1.5). Theorem 1 is proved.

\ \ \ \ \ \ Recall that if $Q$ is a  trigonometric  polynomial on $\T$ of degree at most  $n$,   then
\begin{equation}
\|Q\|_{L^{\infty}(\T)}\le 2n^{1/p}\|Q\|_{L^{p}(\T)}
\end{equation}
for each $1\le p<\infty$  (see [10,  p. 495]). An analogous inequality  holds  in $\be$ (see, for example, [15, p. 233]).  Namely, if $f\in B^{r_1}_{\sigma}$ and  $1\le r_1\le r_2\le\infty$ , then
\begin{equation}
\|f\|_{r_2}\le 2\sigma^{\frac1{r_1}-\frac1{r_2}}\|g\|_{r_1}.
\end{equation}
We will need the following Bernstein inequality
 \begin{equation}
|F(x)-F(y)|\le 2\sin\Bigl(\frac{\sigma|x-y|}{2}\Bigr) \|F\|_{\infty},
\end{equation}
where $F\in\bi$, $x,y\in\R$,  and $\sigma|x-y|\le \pi$ (see [15, p. 213]).

\ \ \ \ \ { PROOF OF THEOREM 2.}\quad Let $f\in \be$, $1< p<\infty$. By (1.5), there exists  a constant  $L(f)>0$ such that
\begin{equation}
\sup _{\tau>0} \|\varphi_{f, \tau}\|_{L^p(-\tau,\tau)}\le L(f)<\infty.
\end{equation}
 Set $u_{f}(t)=f_{\tau}\bigl(\frac{\tau}{\pi}t\bigr)$,  $t\in (-\pi, \pi].$  We will  consider $u_f$ as an element of $L^p(\T)$. Then  (2.15) and (2.18) imply that
\begin{gather}
\|f_{\tau}\|_{L^{\infty}(\R)}=\|f_{\tau}\|_{L^{\infty}([-\tau,\tau])}=\|u_{f}\|_{L^{\infty}(\T)}\le2N^{\frac1{p}}\|u_f\|_{L^p(\T)}=2\biggl(\frac{\pi N}{\tau}\biggr)^{\frac1p}\| f_{\tau}\|_{L^p(-\tau,\tau)}\nonumber \\
\le 2\sigma^{\frac1p}\| f_{\tau}\|_{L^p(-\tau,\tau)}
=2\sigma^{\frac1p}\|\varphi_{f,\tau}\|_{L^p(-\tau,\tau)}\le 2\sigma^{\frac1p}L(f).\nonumber
\end{gather}
 Therefore, using (2.16) in the case where $r_1=p$ and $r_2=\infty$, we have that there is    $M(f)>0$ such that
\begin{equation}
\sup _{\tau>0} \|f- f_{\tau}\|_{L^{\infty}(\R)}\le M(f)<\infty.
\end{equation}
\ \ \ \ \ If $g\in \be$ and $1\le p<\infty$, then $\lim_{t\to\infty}g(t)=0$ (see [7, p. 150]). Hence, in order to prove (1.7),  it suffices  to prove that
$\lim_{\tau\to\infty} \|f-\varphi_{f, \tau}\|_{L^{\infty}(-\tau,\tau)}=0.$
Assume to the contrary that there exist a constant $a>0$ and a sequence $\tau_m\to\infty$ such that
$\|f-\varphi_{f,\tau_m}\|_{L^{\infty}(-\tau_m,\tau_m)}\ge a$,   $m\in \mathbb{N}$. Therefore, we can choose $x_m\in[-\tau_m, \tau_m]$ so that
\begin{equation}
|(f-f_{\tau_m})(x_m)|=|(f-\varphi_{f,\tau_m})(x_m)|\ge a>0,
\end{equation}
$m\in \mathbb{N}$. Now we will estimate $\|f-f_{\tau_m}\|_{L^{p}(-\tau_m,\tau_m)}$. Fix $\delta\in(0,1)$. Assume without loss of generality that $\tau_m\ge \pi/\sigma$ for all  $m\in \mathbb{N}$. Then   for at least one     of the intervals $[x_m-\pi\delta/\sigma, x_m]$ and $[x_m, x_m+\pi\delta/\sigma]$ is contained in  $ [-\tau_m,\tau_m]$.  We denote such interval by $I_m$. For any  $y\in I_m$,  we take this $y$ together with $x=x_m$ and $F=f-f_{\tau_m}$  in (2.17).  Therefore,  if we recall (2.19) and (2.20), we obtain
\begin{equation}
|(f-f_{\tau_m})(y)|\ge |(f-f_{\tau_m})(x_m)| -2 M(f)\sin\frac{\sigma|y-x_m|}{2}\ge a-\pi\delta M(f).
\end{equation}
Also, (2.19) and (2.20) clearly imply $a\le M(f)$.  Thus, substituting $a/(2\pi M(f))$ for $\delta$ in (2.21), we have
\[
\int_{-\tau_m}^{\tau_m}|(f-f_{\tau_m})(y)|^p\,dy\ge \int_{I_m}|(f-f_{\tau_m})(y)|^p\,dy\ge \frac{a^{p+1}}{2^{p+1}\sigma M(f)}
\]
for all $m\in\mathbb{N}$. This contradicts (1.5) and  Theorem 2 is proved.

\ \ \ \  {\bf{ Acknowledgments}}. The author thanks the referee for valuable suggestions which helped to improve this paper.

\newpage
{ \centerline{\bf{References}}
 \vspace{2mm}

[1] \ N. Y. Achieser,\ {\it Theory of Approximation}, Dover Publications. Inc. (New York, 1992).

 [2] \ E. Hewitt and K.A. Ross, \ {\it Abstract Harmonic Analysis II},\ Springer-Verlag (Berlin-Heidelberg, 1997).

[3] \ J.R. Higgins,\ {\it Sampling Theory in Fourier and Signal Analysis: Foundations}, Clarendon Press (Oxford, 1996).

[4] \  L. H\"ormander,\ A new proof and a generalization of an inequality of Bohr, {\it Math. Scand.}, {\bf 2} (1954), 33--45.

[5] \  L. H\"ormander,\ Some inequalities for functions of exponential type, {\it Math. Scand.}, {\bf 3} (1955), 21--27.

[6] \  M.G. Krein,\ On the representation of  functions by Fourier--Stieltjes integrals,   (Russian) {\it Uch. Zap. Kuibysh. Gos. Ped.  Uchit. Inst.}, {\bf 7}  (1943),
 123-147.

[7] \  B.Ya.  Levin,\  {\it Lectures on Entire Functions},\ Transl. Math. Monogr., {\bf 150}, Amer. Math. Soc. (1996).

[8] \  B.M.  Lewitan,\ \"{U}ber eine Verallgemeinerung der Ungleichungen von    S. Bernstein und  H. Bohr,   {\it  Dokl. Akad. Nauk}., {\bf 15}  (1937), 169--172.

[9] \ R. Martin,\ Approximation of\ $\Omega$--bandlimited functions by\ $\Omega$--bandlimited trigonometric polynomials, {\it Sampl. Theory Signal Image Process}., {\bf 6}(3) (2007), 273--296.

[10] \ G.V. Milovanovi\'{c},   D.S. Mitrinovi\'{c} and  Th. M. Rassias,  {\it Topics in Polynomials: Extremal problems, inequalities, Zeros},  World Scientific Publishing Co., Inc.
(River Edge, NJ, 1994).

[11] \ S. Norvidas,\ Approximation of entire functions by exponential polynomials, {\it Lith. Math. J}., {\bf 34}(4)  (1994), 415--421.

[12] \ S. Norvidas,\ Approximation of bandlimited  functions by finite exponential sums, {\it Lith. Math. J}., {\bf 49}(2) (2009), 185--189.

 [13] \  Q.I. Rahman and G. Schmeisser,\  $L^p$ inequalities for entire functions of exponential type, {\it Trans. Amer. Math. Soc}., {\bf 320} (1990), 91--103.

[14] \ G. Schmeisser,\ Approximation of entire functions of exponential type by trigonometric polynomials,  {\it Sampl. Theory Signal Image
Process}., {\bf 6}(3) (2007), 297--306.

[15] \ A.F. Timan,\ {\it Theory of Approximation of Functions of a Real Variable},  Dover Publications, Inc. (New York, 1994).

[16] \ R.M. Trigub and E.S. Bellinsky,\ {\it Fourier Analysis and Approximation of Functions}, Kluwer Academic Publishers  (Dordrecht, 2004).
}}
\end{document}